\documentstyle[a4,amsfonts,12pt]{article}

\newcommand{\mb}{\mbox}
\newcommand{\beq}{\begin{equation}}
\newcommand{\eeq}{\end{equation}}
\newcommand{\ueberschrift}{\bigskip\goodbreak\noindent\bigskip}
\newcounter{theabsatz}
\newcommand{\absatz}[1]{\stepcounter{theabsatz} \ueberschrift
               {\large \bf \arabic{theabsatz}. {#1}} \setcounter{equation}{0}}

 \newtheorem{lem}{Lemma}
 
\newcommand{\z}{\zeta }
\begin{document}
\newcommand{\ext}{\mbox{ext\,}}
\newcommand{\diam}{\mbox{diam\,}}

\parindent 10 pt
\parskip 8pt plus 4pt
\jot 10pt

\abovedisplayskip 8pt plus 1pt \belowdisplayskip 8pt plus 1pt

\newcommand{\C}{{\mbox{C}}}
\newcommand{\OC}{\overline{\mbox{\bf C}}}
\newcommand{\Z}{{\cal Z}}
\newcommand{\CC}{{\overline{\mbox{C}}}}
\newcommand{\DD}{{\overline{\mbox{D}}}}
\newcommand{\D}{{\mbox{D}}}
\newcommand{\De}{\Delta}
\newcommand{\R}{{\mbox{R}}}
\newcommand{\T}{{\mbox{T}}}
\newcommand{\N}{{\mbox{N}}}
\newcommand{\PP}{{\mbox{P}}}
\newcommand{\He}{{\mbox{H}}}
\newcommand{\p}{\preceq}
\newcommand{\s}{\succeq}
\newcommand{\En}{{\mathcal E}_n}
\newcommand{\ov}{\overline}
\newcommand{\de}{\delta}
\newcommand{\Ga}{\Gamma}
\newcommand{\ga}{\gamma}
\newcommand{\La}{\Lambda}
\newcommand{\oge}{\succeq}
\newcommand{\ole}{\preceq}
\newcommand{\la}{\lambda}
\newcommand{\kap}{\mb{cap }}
\newcommand{\be}{\beta}
\newcommand{\Om}{\Omega}
\newcommand{\Si}{\Sigma}
\newcommand{\om}{\omega}
\newcommand{\bP}{{{\mathbb{P}}}}
%{{\bf \Pi}}
\newcommand{\al}{\alpha}
\newcommand{\ve}{\varepsilon}

\begin{center}
{\large \bf   On Lebesgue Constants for Interpolation Points 
on a Quasiconformal Arc}\\[3ex]

 {\bf Vladimir Andrievskii}\\[3ex]

{\it  Department of Mathematical Sciences, Kent State University,
Kent, OH 44242}
 \end{center}

\vspace*{2cm}

Running head: Lebesgue constants

\vspace*{2cm}

Mailing address:

\noindent V.V. Andrievskii\\   Department of Mathematical Sciences\\
Kent
 State University \\ Kent, OH 44242, USA\\[4ex]

{\it E-mail address}:\\ andriyev@math.kent.edu\\

{\it Phone}: (330) 672 9029

 \newpage

\begin{center}
 {\bf Abstract}
 \end{center}

 Using the theory of quasiconformal mappings, we simplify the proof of the
recent result by Taylor and Totik (see IMA Journal of Numerical Analysis
{\bf 30} (2010) 462--486) on the behavior of the Lebesgue constants
for interpolation points
on a compact set in the complex plane.

\vspace*{0.5cm}



 \newpage

\absatz{Introduction}

For  a compact set $K$ in the complex plane $\C$ with a positive logarithmic capacity and
$\N:=\{1,2,\ldots\}$, consider a triangular array $Z:=\{z_{n,k}\}_{1\le k\le n,n\in\N}$
of points in $K$ (i.e., an {\it interpolation scheme} on $K$) such that
$z_{n,k}\neq z_{n,j}$ for $k\neq j$.
In the theory of the Lagrange interpolation of continuous functions on $K$ an important role is played by 
the 
{\it Lebesgue constants}
$$
\Lambda_n=\Lambda_n(K,Z):=\sup_{z\in K}\sum_{k=1}^n|l_{n,k}(z)|,
$$
where
$$
l_{n,k}(z):=\frac{\prod_{1\le j\le n,j\neq k}(z-z_{n,j})}{\prod_{1\le j\le n,j\neq k}(z_{n,k}-z_{n,j})}
$$
(see \cite{gai, smileb}).

The starting point of our consideration are the following recent results by Taylor
and Totik \cite{taytot}. 
We refer the reader to \cite{ran, saftot} for basic notions of potential theory.
We say that $\La_n$ are {\it subexponential} if
\beq\label{1.1}
\La_n^{1/n}\to 1,\quad \mb{for }n\to\infty.
\eeq
We associate with the rows of $Z$ the {\it normalized counting measures}
$$
\nu_n:=\frac{1}{n}\sum_{k=1}^n\de_{z_{n,k}},\quad n\in\N,
$$
where $\de_z$ denotes the unit mass at $z\in\C$.

Let $\mu_K$ be the equilibrium measure of $K$. We say that
$\nu_n$ are {\it asymptotically distributed like }$\mu_K$ if
$\nu_n\to \mu_K$ in the weak$^*$ sense.

{\bf Theorem} A (\cite[Theorem 1.1]{taytot}). {\it If $\La_n$ are subexponential then
$\nu_n$ are asymptotically distributed like $\mu_K$}.

Simple examples show that even for $K=[-1,1]$, the inverse 
is not necessarily true.

{\bf Theorem} B (\cite[Theorem 6.3]{taytot}). {\it Let $K=[-1,1]$. If 
$\nu_n$ are asymptotically distributed like $\mu_K$ and if each pair of points in the
nth row of $Z$ satisfies the distancing rule
\beq\label{1.2}
|z_{n,j}-z_{n,k}|\ge c_1\left(\frac{\sqrt{1-|z_{n,j}|}}{n}+\frac{\sqrt{1-|z_{n,k}|}}{n}+
\frac{1}{n^2}\right),\quad 1\le j,k\le n,j\neq k
\eeq
for some constant $c_1>0$, then $\La_n$ are subexponential.}

The last theorem plays a crucial role in establishing (\ref{1.1}) for the Leja points on
the compact set with a piecewise $C^2$-smooth outer boundary. 

The main objective of our paper is to extend Theorem B to the case of an arbitrary quasiconformal arc
(see \cite{ahl, lehvir})
and simplify
 its original proof.
We use the following traditional for the approximation theory in the complex plane idea.

For  $K=I:=[-1,1]$
denote by $I_{1/n}, n\in\N$, the ellipse with foci at $\pm 1$
and the sum of semiaxes equal to $1+1/n$. Such an ellipse is the image
of the circle $\{ w:\, |w|=1+1/n\}$ under the Joukowski mapping
$z=\Psi(w)=(w+1/w)/2$ of $\D^*:=\{ w:\, |w|>1\}$
onto $\OC\setminus I$, where $\OC:=\C\cup\{\infty\}$. 
Then, for $x\in I$ and $n\in\N$,
$$
\frac{1}{c_2}\rho_{1/n}(x)\le\frac{\sqrt{1-|x|}}{n}+\frac{1}{n^2}\le c_2\rho_{1/n}(x)
$$
holds with a constant $c_2\ge 1$, where
$$
\rho_{1/n}(x):= \inf_{\z\in I_{1/n}}|\z-x|,\quad x\in I
$$
is the distance from $x$ to $I_{1/n}$

The notions of $\Psi,I_{1/n}$ and $\rho_{1/n}$ are
 also meaningful for a 
 bounded arc $K=L$ in $\C$ which is the key to a
generalization
of Theorem B. 

\absatz{Main result}

From now on, we make the assumption  that $\Om:=\CC\setminus K$ is connected and regular with respect
to the Dirichlet problem. Denote by $g_\Om$ the Green function of $\Om $ with pole at
$\infty$, and let for $z\in K$ and $\de>0$,
$$
K_\de:=\{\z\in\Om:g_\Om(\z)=\log(1+\de)\},\quad
\rho_\de(z):=d(z,K_\de),
$$
where
$$
d(z,B):=\inf_{\z\in B}|z-\z|,\quad z\in \C,B\subset\C.
$$
We say that $Z$ is {\it well separated} (on $\partial K$) if
\beq\label{2.0}
z_{n,k}\in \partial K,\quad n,k\in\N,1\le k\le n
\eeq
and each pair of points in the $n$th row satisfies 
\beq\label{2.1}
|z_{n,k}-z_{n,j}|\ge c_3\rho_{1/n}(z_{n,k}),\quad 1\le k,j\le n,k\neq j,
\eeq
where $c_3>0$ is a constant.

In numerical analysis one of the most important  examples of $Z$ with (\ref{2.0})-(\ref{2.1})
is the set of {\it Leja points} which is defined inductively as follows.
Let $z_1\in \partial K$ be arbitrary. If $z_1,\ldots,z_{n-1}$ are known,
then let $z_n\in \partial K$ is a point at which 
$\prod_{j=1}^{n-1}|z-z_j|$ attains its maximum. The row $\{z_{n,k}\}_{k=1}^n$
consists of the first $n$ Leja points.
Note that $\nu_n$ for the Leja points
are asymptotically distributed like $\mu_K$ (see \cite[Chapter V.1]{saftot}).

{\bf Proposition.} {\it Leja points are well separated}.

In what follows, $L$ denotes a bounded {\it quasiconformal} arc which means that for any 
two points
$z,\z\in L$,
\beq\label{2.2}
\mb{diam }L(z,\z)\le c_4|z-\z|,
\eeq
where $c_4\ge 1$ is a constant, $L(z,\z)$ is a subrarc of $L$ between these points, and
 diam $S$ is the diameter of a set
$S\subset\C$.

{\bf Theorem.} {\it Let $K=L$. If $\nu_n$ are asymptotically distributed like $\mu_L$ and if
$Z$ is well separated, then $\La_n$ are subexponential.}

In the case $L=[-1,1]$, the above theorem implies Theorem B.
Using the reasoning of \cite[Section 7]{taytot}, one can further extend the
above theorem to the case of $K$ whose boundary consists of a finite
number of quasiconformal arcs. 
Reichel
(see \cite[p. 465]{taytot}), based
on a numerical calculation, conjectured 
that the Lebesgue constants for Leja points are
subexponential also for an arbitrary $K$.
This challenging conjecture remains open still.

We denote by $c,c_1,c_2,\ldots$ positive constants (different in different sections)
which can depend only on $K,L,$ and $Z$. Moreover, for positive functions $a$ and $b$ we use
the order inequality $a\ole b$ if $a\le cb $. The expression $a\asymp b$ means
that $a\ole b$  and $b\ole a$ simultaneously.

\absatz{Auxiliary results}

Let $K=L$ and let function $\Phi$ map $\Om$ conformally and univalently  onto $\D^*$ such that
$\Phi(\infty)=\infty, \Phi'(\infty)>0$. Set $\Psi:=\Phi^{-1}$. Since in this case
$g_\Om=\log|\Phi|$, we have
$$
L_\de=\{\z\in\Om:|\Phi(\z)|=1+\de\},\quad \de>0.
$$

In this section, we mention some known results about metric properties
of $\Phi$ and $\Psi$ (see for more details \cite{andbla}) and their consequences.
Denote by $z_1$ and $z_2$ the endpoints of $L$. Since $\Phi$ can be extended continuously to
these points, we set for $\de>0, j=1,2$, and $z\in L$
$$
t_j:=\Phi(z_j),\quad \D_1^*:=\{t:|t|>1, \arg t_1<\arg t<\arg t_2\},
$$
$$
\D_2^*:=\D^*\setminus \ov{\D_1^*},\quad \Om_j:=\Psi(\D^*_j),
$$
$$ 
L^j_\de:=L_\de\cap\ov{\Om_j},\quad \rho_\de^j(z):=d(z,L_\de^j),
$$
so that
$$
\rho_\de(z)=\min_{j=1,2}\rho_\de^j(z).
$$
Denote by $\Phi_j,j=1,2$, the restriction of $\Phi$ to $\Om_j $
and let $\tilde{z}^j_\de:=\Psi((1+\de)\Phi_j(z))$.
Since $\partial\D^*_j$ and $\partial\Om_j$ are quasiconformal
(see \cite[p. 30, Lemma 2.8]{andbla}),
$\Phi_j$ can be extended to a quasiconformal mapping $\Phi_j:\OC\to\OC$.
Therefore, the following statement holds.
\begin{lem}\label{lem3.1} (\cite[p. 29, Theorem 2.7]{andbla})
Let $\z_k\in\ov{\Om_j}\setminus\{\infty\},\Phi_j(\z_k)=:\tau_{j,k},j=1,2,k=1,2,3.$
Then the conditions $|\z_1-\z_2|\ole |\z_1-\z_3|$ and 
$|\tau_{j,1}-\tau_{j,2}|\ole |\tau_{j,1}-\tau_{j,3}|$
are equivalent.
\end{lem}
Let $z_\de^j\in L_\de^j$ satisfy $\rho_\de^j(z)=|z-z_\de^j|$.
Lemma \ref{lem3.1} with the triplet of points $z,\tilde{z}^j_\de, z_\de^j$
implies
\beq\label{3.1}
\rho_\de^j(z)\asymp|z-\tilde{z}^j_\de|.
\eeq
Let $\T:=\{t:|t|=1\}$ be the unit circle and $|S|$ denote the length of an arc $S\subset \C$.
The function $\Psi$ can be extended
 continuously to a function $\Psi:\ov{\D^*}\to\OC$, and for any subarc $J\subset L$, there exist two arcs
 $J_1'\subset\ov{\D^*_1}\cap\T$ and $J_2'\subset\ov{\D^*_2}\cap\T$ such that
 $\Psi(J_1')=\Psi(J_2')=J$ and  the intersection $J_1'\cap J_2'$ consists of at most two points.
 Then
 \beq\label{3.2}
 \mu_L(J)=\frac{1}{2\pi}\left(|J_1'|+|J_2'|\right)
 \eeq
(see \cite[p. 22]{andbla}).
\begin{lem} \label{lem3.2}
 For $Z=Z(L)$ the spacing condition (\ref{2.1}) implies
$$
\mu_L(L(z_{n,k},z_{n,j}))\oge\frac{1}{n},\quad
1\le k,j\le n,k\neq j.
$$
\end{lem}
{\bf Proof.} Let $j_0=j_0(z_{n,k},n)$ satisfy 
$\rho_{1/n}(z_{n,k})=\rho_{1/n}^{j_0}(z_{n,k}).$
Since by (\ref{2.1}) and (\ref{3.1})
$$
|z_{n,k}-\tilde{(z_{n,k})}^{j_0}_{1/n}|\asymp
\rho_{1/n}^{j_0}(z_{n,k})=\rho_{1/n}(z_{n,k})
\ole|z_{n,k}-z_{n,j}|,
$$
Lemma \ref{lem3.1} with the triplet of points $z_{n,k}, \tilde{(z_{n,k})}^{j_0}_{1/n},
z_{n,j}$ and (\ref{3.2}) yield
$$
\mu_L(L(z_{n,k},z_{n,j}))\ge \frac{1}{2\pi}|\Phi_{j_0}(z_{n,k})-\Phi_{j_0}(z_{n,j})|
\oge\frac{1}{n}.
$$

\hfill$\Box$

In the proof of the next lemma we use the notion and properties of the
module of a family of curves, see \cite{ahl, lehvir, andbla} for more details.

\begin{lem}\label{lem3.3}
For $\xi_1,\xi_2\in L$,
\beq\label{3.4}
\mu_L(L(\xi_1,\xi_2))\ole |\xi_1-\xi_2|^{1/2}.
\eeq
\end{lem}
{\bf Proof.} Let
$J:=L(\xi_1,\xi_2).$
Without loss of generality we can assume that \newline
diam $J=: d< D:=($diam $L)/2$.

The main idea of the proof is to compare modules of the following families of curves.
Denote by $\Ga_j',j=1,2$ the family of all crosscuts of $\D^*$ which separate 
$J_j':=\Phi_j(J)$ from $\infty$. Let $\Ga_j:=\Psi(\Ga_j')$ and 
$$
\Ga_3:=\left\{\ga_r:=\{\z:|\z-\xi_1|=r\}:d<r< D\right\}.
$$
Since $\Ga_j<\Ga_3$, i.e., each $\ga_3\in\Ga_3$ contains $\ga_j\in\Ga_j$,
according to \cite[p. 343, Theorem 1.2; pp. 347-349, Examples 1.9 and 1.11]{andbla},
for their modules we have
$$
\frac{1}{2\pi}\log\frac{D}{d}=m(\Ga_3)\le m(\Ga_j)=m(\Ga_j')
\le 2+\frac{1}{\pi}\log\frac{4}{|J_j'|}.
$$
Therefore, (\ref{2.2}) yields
$$
|J_j'|\le \frac{4 e^{2\pi}}{D^{1/2}}d^{1/2}\ole |\xi_1-\xi_2|^{1/2},
$$
which, by virtue of (\ref{3.2}), implies (\ref{3.4}).

\hfill$\Box$

\absatz{Proofs}

{\bf Proof of Proposition.}
We use mathematical induction. Let $\{z_j\}_{j\in\N}$ be the
Leja points for $K$. For $n=2$ (\ref{2.1}) is trivial. Next, assuming that 
(\ref{2.1}) is true for $m=n-1$, consider the polynomial
$$
p_{n-1}(z):=\prod_{j=1}^{n-1}(z-z_j)
$$
for which we have $||p_{n-1}||_K=|p_{n-1}(z_n)|.$
Here $||\cdot||_S$ is the uniform norm on $S\subset\C$.
By the Bernstein-Walsh Lemma (see \cite[p. 77]{wal} or \cite[p. 153]{saftot})
$$
||p_{n-1}||_{K_{1/n}}\le\left(1+\frac{1}{n}\right)^{n-1}||p_{n-1}||_K
<e||p_{n-1}||_K.
$$
Therefore, for $z\in\partial K$ and $\z\in\C$ with $|\z-z|\le \rho/2$,
where $\rho:=\rho_{1/n}(z)$, we have
$$
|p_{n-1}'(\z)|\le\frac{1}{2\pi}\int_{\{\xi:|\xi-\z|=\rho/2\}}
\frac{|p_{n-1}(\xi)|}{|\xi-\z|^2}|d\xi|
\le\frac{2e}{\rho}||p_{n-1}||_K.
$$
To verify (\ref{2.1}) for $m=n$ we can assume that
one of points in the left-hand side of (\ref{2.1}) is $z_n$.
In the case where $|z_n-z_j|<\rho_{1/n}(z_n)/2,j<n$
we obtain 
$$
||p_{n-1}||_K=|p_{n-1}(z_n)|\le\int_{[z_n,z_j]}|p'_{n-1}(\z)||d\z|
\le \frac{2e|z_n-z_j|}{\rho_{1/n}(z_n)}||p_{n-1}||_K
$$
from which we have $|z_n-z_j|\oge\rho_{1/n}(z_n)$.
Moreover,
$$
\rho_{1/n}(z_j)\le|z_j-z_n|+\rho_{1/n}(z_n)\ole |z_j-z_n|
$$
which completes the proof of 
(\ref{2.1}).

\hfill$\Box$

{\bf Proof of Theorem.}
From the reasoning in \cite[Section 4]{taytot} we see that it is enough to show that for
any sufficiently small $\de>0$ there exists $c_\de>0$ such that
\beq\label{4.1}
\lim_{\de\to0^+}c_\de=0,
\eeq
\beq\label{4.2}
\liminf_{n\to\infty}\min_{1\le k\le n}S_{n,k,\de}^{1/n}\ge e^{-c_\de},
\eeq
where
$$
S_{n,k,\de}:=\prod_{z_{n,j}\in A_{n,k,\de}}
|z_{n,j}-z_{n,k}|
$$
and
$$
A_{n,k,\de}:=\{z_{n,j}:1\le j\le n,j\neq k,|z_{n,j}-z_{n,k}|\le\de\}.
$$
Let $\z_0:= z_{n,k}$. In the most complicated case where 
$A_{n,k,\de}\cap L(\z_0,z_j)\neq\emptyset, j=1,2$ (here, as before,
$z_1$ and $z_2$ are the endpoints of $L$), which we consider in detail,
we rename the points in $A_{n,k,\de}$ as follows
$$
A_{n,k,\de}=\bigcup_{m=-m_1,m\neq 0}^{m_2}\{\z_m\},
$$
where 
$\z_m\in L(\z_{m-1},\z_{m+1}), -m_1+1\le m\le m_2-1$.

Since by Lemma \ref{lem3.3}
$$
\mu_L(L(\z_{-m_1},\z_{m_2}))\ole\de^{1/2},
$$
Lemma \ref{lem3.2} implies
$$
m_j\le c_1\de^{1/2}n,\quad j=1,2.
$$
Taking into account that
$$
m!\ge\left(\frac{m}{e}\right)^m,\quad m\in\N,
$$
and by Lemmas \ref{lem3.2} and \ref{lem3.3} for $-m_1\le m\le m_2, m\neq 0$,
$$
|\z_m-\z_0|
\ge c_2\mu_L(L(\z_m,\z_0))^2\ge \left( \frac{c_3| m|}{n}\right)^2
$$
we have
$$
S_{n,k,\de}\ge
\prod_{m=-m_1,m\neq 0}^{m_2}\left(\frac{c_3| m|}{n}\right)^2
\ge \left(\frac{c_3m_1}{en}\right)^{2m_1}
\left(\frac{c_3m_2}{en}\right)^{2m_2}.
$$
Therefore, for sufficiently small $\de$,
$$
S_{n,k,\de}^{1/n}\ge
\left(\frac{c_1c_3\de^{1/2}}{e}\right)^{4c_1\de^{1/2}}
$$
which implies (\ref{4.2}) with
$$
c_\de:=4c_1\de^{1/2}\log\frac{e}{c_1c_3\de^{1/2}}
$$
satisfying (\ref{4.1}).

\hfill$\Box$

\absatz{Acknowledgements}

Part of this work was done during the Fall of 2016 semester, while the author visited
the Katholische Universit\"at Eichst\"att-Ingolstadt and
the Julius Maximilian University of W\"urzburg.
The author is  also grateful to   M.
 Nesterenko
 for his helpful comments.

\end{document}